\documentclass[12pt]{amsart}
\usepackage{amscd}
\usepackage{amsmath}
\usepackage{amssymb}
\usepackage{slashbox}
\usepackage{hyperref}
\usepackage{graphicx}
\usepackage{xcolor}
\usepackage{tikz}

\def\Null{\operatorname{Null}}

\def\supp{\operatorname{supp}}

\newcommand{\ZZ}{\mathbb Z}
\newcommand{\NN}{\mathbb N}

\newcommand{\C}{\mathcal C}
\newcommand{\QQ}{\mathbb Q}

\newtheorem{lemma}{Lemma}[section]
\newtheorem{corollary}[lemma]{Corollary}
\newtheorem{theorem}[lemma]{Theorem}
\newtheorem{proposition}[lemma]{Proposition}
\newtheorem{definition}[lemma]{Definition}
\newtheorem{remark}[lemma]{Remark}
\newtheorem{example}[lemma]{Example}

\newtheorem{setup}[lemma]{Setup}

\advance\headheight1.15pt

\begin{document}

\title[Strongly robust toric ideals of weighted oriented graphs]{Strongly robustness of toric ideals of weighted oriented graphs}

\author[R. Nanduri]{Ramakrishna Nanduri$^*$}
\address{Department of Mathematics, Indian Institute of Technology
Kharagpur, West Bengal, INDIA - 721302.}
\email{nanduri@maths.iitkgp.ac.in}
\author[T. K. Roy]{Tapas Kumar Roy$^{\dag}$}
\address{Department of Mathematics, Indian Institute of Technology
Kharagpur, West Bengal, INDIA - 721302.}
\email{tapasroy147@kgpian.iitkgp.ac.in, tapasroy147@gmail.com}
\thanks{$^\dag$ Supported by PMRF research fellowship, India.}
\thanks{$^*$ Corresponding author}
\thanks{{\bf AMS Classification 2020:} 13F65, 05E40, 05C22, 20M25, 13A70.}
\maketitle

\begin{abstract}
 In this article, we investigate the strongly robust property of toric ideals associated with weighted oriented graphs. We establish that the toric ideals of a broad class of monomial ideals are strongly robust; this class encompasses the edge ideals of weighted oriented graphs in which every edge is incident to a vertex of degree $2$.
\end{abstract}

\section{Introduction} \label{sec1}

Toric ideals play a central role in the modern theory of Commutative Algebra, Combinatorial Commutative Algebra, Combinatorial Optimization, and Algebraic Coding Theory and related fields (see \cite{s95}, \cite{v01}). Consider an $ n\times m$ matrix $A=[{\bf a_1} \ldots {\bf a_m}]$ with columns ${\bf a_1}, \dots, {\bf a_m}$ consisting of non-negative integer entries. Let $M=({\bf x^{a_1}, \ldots, {\bf x^{a_m}}})$ be the be the corresponding monomial ideal in the polynomial ring $R=K[x_1,\ldots,x_n]$, where $K$ is a field. Let $S=K[e_1,\ldots,e_m]$ be a polynomial ring in the variables $e_1,\dots,e_m$. Define a ring homomorphism $\varphi: S \rightarrow R$, by $\varphi(e_i)= {\bf x^{a_i}}$. The kernel of $\varphi$, denoted $I_A$, is called the {\it toric ideal} associated to $A$ or $M$. Note that $S/I_A$ is the affine semigroup ring corresponding to the semigroup $\ZZ_+{\bf a_1}+\cdots +\ZZ_+{\bf a_m}$. It is well known that $I_A$ is a prime ideal generated by binomials ${\bf e^u}-{\bf e^v}$ such that $\deg_A({\bf e^u})=\deg_A({\bf e^v})$, where $\deg_A({\bf e^u}) =u_1{\bf a_1}+\cdots+u_n{\bf a_n}$.  

A binomial ${\bf e^u}-{\bf e^v} \in I_A$ is called a {\it primitive} binomial if there exists no other binomial ${\bf e^{u^{\prime}}}-{\bf e^{v^{\prime}}} \in I_A$ such that ${\bf e^{u^{\prime}}}\vert {\bf e^u}$ and ${\bf e^{v^{\prime}}} \vert {\bf e^v}$. The set of all primitive binomials in $I_A$, denoted by $Gr_A$, is called the {\it Graver basis} of $I_A$. Note that $Gr_A$ generates $I_A$, but it is generally not a minimal generating set. The minimal binomial generating systems of toric ideals have been extensively studied (see \cite{s95}, \cite{ckt07}). In algebraic statistics, a key question is whether a toric ideal has a unique minimal system of binomial generators (see \cite{ta04}). To investigate this, researchers have introduced various notions for binomial generators of toric ideals, including robust, strongly robust, and weakly robust properties. A toric ideal $I_A$ is {\it robust} if it is minimally generated by its universal Gr\"obner basis $\mathcal{U}_A$, defined as the union of all reduced Gr\"obner bases of $I_A$ with respect to all possible term orders. Boocher and Robeva \cite{br15} introduced robust toric ideals and proved that those generated by quadrics are essentially determinantal. Robust toric ideals exhibit many interesting properties; for instance, their multigraded Betti numbers coincide with those of any initial ideal $\text{in}_<(I_A)$, for every term order $<$. Later Sullivant \cite{s19} defined a stronger notion of robustness.  

A toric ideal $I_A$ is {\it strongly robust} if its Graver basis is a minimal generating set of $I_A$, equivalently, if $Gr_A$ equals the set of all indispensable binomials of $I_A$. Note that every strongly robust toric ideal is robust. A binomial $f={\bf e^u}-{\bf e^v}$ in $I_A$ is {\it indispensable} if every system of binomial generators of $I_A$ contains $f$ or $-f$. Surprisingly, robustness and strong robustness coincide for the toric ideals of simple graphs \cite{bbdlmns15}; see also \cite{t16}. Studies on strongly robust toric ideals appear in \cite{oh99, oh05, v17, s19, nn22, mot22, ktv23}. Nevertheless, the notions of robust and strongly robust properties remain completely unexplored for toric ideals of weighted oriented graphs. 

In this work, we characterize the strongly robust property of the toric ideal associated to the edge ideal of a vertex-weighted oriented graph (which we call the toric ideal of the weighted oriented graph). A weighted oriented graph $D$, denoted shortly as a WOG, is a triple $D=(V(D), E(D), w)$, where $V(D)=\{x_1,\ldots,x_n\}$ is the finite set of vertices of $D$, $E(D)$ is the set of edges of $D$ (with no multiple edges or loops), defined as  
$$E(D)=\{(x_i,x_j): \text{there is a directed edge from $x_i$ to $x_j$} \},$$ 
and $w:V(D) \rightarrow \NN$, is a weight function assigning a positive integer weight to each vertex of $D$. Each edge in $E(D)$ carries a specified orientation. 

Weighted oriented graphs (WOGs) are significant due to their broad applications across various research fields; see \cite{mpv17, gmsvv18, hlmrv19, prt19, zxwt19, kblo22, kn23} for a thorough study. The edge ideal of a weighted oriented graph $D$ on $n$ vertices $\{x_1,\ldots,x_n\}$ is defined as the monomial ideal  
$$I(D)= \left(x_ix_j^{w_j} : (x_i, x_j)\in E(D) \right) \subset R.$$ 
The toric ideal corresponding to the edge ideal $I(D)$ of $D$ is denoted by $I_D$. This toric ideal $I_D$ was studied in \cite{bklo}, where it was characterized when principal. Let $\mathcal{M}_D, \mathcal{U}_D $ and $Gr_D$ denote the universal Markov basis, the universal Gr\"obner basis and the Graver basis of $I_D$, respectively. An irreducible binomial ${\bf e^u}-{\bf e^v} \in I_A$ is a {\it circuit} of $I_D$ if it has minimal support with respect to set-inclusion. The set of all circuits of $I_D$ is denoted by $\mathcal{C}_D$. We have  
$\mathcal{C}_D \subseteq \mathcal{U}_D  \subseteq Gr_D \text{ and } \mathcal{M}_D \subseteq Gr_D.$

In \cite{nr23}, the authors completely characterized the circuit binomials of any WOG and explicitly computed them in terms of the minors of its incidence matrix. In this work, we prove that a certain class of monomial ideals has strongly robust toric ideals. 
Let $M$ be a monomial ideal in $R$ with the minimal generating set ${\mathcal{G}}(M)=\{{\bf x}^{{\bf a}_{i}} : 1\le i\le m\}$. Assume that $|\supp({\bf x}^{{\bf a}_{i}})|=2$ for each $i=1,2,\ldots,m$. Suppose further that for every ${\bf x}^{{\bf a}_{i}}\in {\mathcal{G}}(M)$, there exists a variable $x_l$ in $\supp({\bf x}^{{\bf a}_{i}})$ and another monomial ${\bf x}^{{\bf a}_{j}} \not = {\bf x}^{{\bf a}_{i}}$ in $\mathcal{G}(M)$ such that $x_l \in {\supp}({\bf x}^{{\bf a}_{i}})\cap {\supp}({\bf x}^{{\bf a}_{j}})$ and $x_l \notin {\supp}({\bf x}^{{\bf a}_{k}})$, for all $k\neq i,j$. Then the toric ideal $I_{M}$ is strongly robust (see Theorem \ref{sec3thm2}). As a consequence, we show that if $D$ is a WOG such that every edge meets a degree $2$ vertex, then $I_D$ is strongly robust (Corollary \ref{sec3cor4}). 
\vskip 0.2cm 
\noindent 
Let $D$ be a weighted oriented graph with no cycles sharing a path, where each edge of every cycle meets a degree 2 vertex. Assume that for any cycles ${\C_1}, {\C_2}, {\C_3}$ in $D$ such that ${\C_1}$ and ${\C_2}$ are connected to ${\C_3}$ by paths $P_1$ and $P_2$, respectively, every edge $e\in E(P_1)\cap E(P_2)$ satisfies one of the following conditions:\\
$(i)$\; $e$ meets a degree 2 vertex; \\
$(ii)$\; ${\C_3}$ shares no vertex with any other cycle in $D$ and is not connected by a path $P$ to any cycle in $D$ with $e\notin E(P)$. \\
Then $I_{D^{}}$ is strongly robust (Theorem \ref{sec3thm4}). As an application, consider the following class of weighted oriented graphs (WOGs) whose toric ideals are strongly robust: 
let $D$ consist of weighted oriented cycles sharing a common vertex $v$. Let $D^{\prime}$ be the WOG obtained by adding new cycles ${\C_1},{\C_2},\cdots,{\C_k}$ disjoint from $D$, where each ${\C_i}$ connects to $v$ via a path $P_i$ with $E(P_i)\cap E(P_j)=\emptyset$ for all distinct $i,j\in\{1,2,\ldots,k\}$. Then $I_{D^{\prime}}$ is strongly robust (Corollary \ref{cor10}).  
 
\vskip 0.2cm
We organize the paper as follows. Section \ref{sec2} recalls the necessary definitions and preliminary results for our main theorems. Section \ref{sec3} then establishes these main results.   

\section{Preliminaries} \label{sec2}

Let $R=K[x_1,\ldots,x_n]$, and $S=K[e_1, \ldots ,e_m]$ be polynomial rings, where $K$ is a field. Recall that the {\it support} of a binomial $ f=\displaystyle \prod_{k=1}^{m} e_k^{p_k} - \prod_{k=1}^{m} e_k^{q_k}$, denoted as $\supp(f)$, is the set of variables $e_k$ which appear in $f$ with non-zero exponent. 

Let $D= (V(D),E(D),{\bf w})$, be a weighted oriented graph. Let $V(D)=\{x_1, \ldots,x_n\}$ and $E(D)=\{e_1,\ldots,e_m\}$. For simplicity, we set $w_j:=w(x_j)$ for each $j$. 

 \begin{definition} \label{def1}
 The {\it incidence matrix} of $D$ is an $n \times m$ matrix whose $(i,j)^{th}$ entry $a_{i,j}$ is defined by
\begin{center}
$a_{i,j}=
\begin{cases}
1\; \mbox { if } \;e_j = (x_i, x_l) \in E(D)\; \mbox{ for some }\;  1 \le l \le n,  \\
w_i\; \mbox{ if }\;e_j = (x_l, x_i) \in E(D)\; \mbox{ for some }\; 1 \le l \le n, \\ 0\; \mbox{ otherwise, }
\end{cases}
$
\end{center}
and we denote by $A(D)$. 
\end{definition}

A weighted oriented $r$-cycle denoted by $\C_r$ is a weighted oriented graph whose underlying graph is a cycle of length $r$ and the edges oriented in any direction. The degree of a vertex $v$ is the number of edges incident (inward or outward) to $v$. A vertex of degree $1$ is called a leaf. Recall that a weighted oriented even cycle $\mathcal{C}_m$ on $m$ vertices is said to be {\it balanced} if $\displaystyle \prod_{k=1}^{m} a_{k,k} =a_{1,m} \prod_{k=1}^{m-1} a_{k+1,k}$, that is, $\mbox{det}(A({\C_m}))=0$, where $A(\mathcal{C}_m)=[a_{i,j}]_{m \times m}$. We denote Null$(A(D))$, the null space of $A(D)$ over $\QQ$. 

  
\begin{definition}
For a vector $\mathbf{b} = \left(b_1,b_2, \ldots,b_{l}\right) \in \ZZ^{l}$, define the corresponding binomial in the variables $e_1,\ldots,e_l$ as
$f_{\mathbf{b}} :=f_{\mathbf{b}}^{+} -f_{\mathbf{b}}^{-},$  where
$$f_{\mathbf{b}}^{+} = \prod\limits_{b_{i}\geq 0}^{l}e_i^{b_i}, \mbox{ and } f_{\mathbf{b}}^{-} = \prod\limits_{b_{i}<0}^{l}e_i^{-b_i}.$$
\end{definition}
\noindent 
For any vector ${\bf b} \in \ZZ^{l}$, we can always write ${\bf b}={\bf b}_{+}-{\bf b}_{-}$ uniquely, for some ${\bf b}_{+}, {\bf b}_{-} \in \ZZ_+^{l}$. Let $[{\bf b}]_{i}$ denotes the $i$-th entry of ${\bf b}$. Also, the support of a vector ${\bf b}$ is defined as supp$({\bf b}) := \{i : [{\bf b}]_{i}\neq 0\}$. If $f_{\bf m}\in I_D$, then ${\bf m}\in \text{Null}(A(D))$ and $[{\bf m}]_{i}$ denotes the $i$-th entry of ${\bf m}$ corresponding to the edge $e_i\in E(D)$.  

\begin{remark} \label{sec2rmk1} 
Let ${\bf b} \in \ZZ^l$. Then it is easy to see that 
\begin{enumerate}
    \item if $[{\bf b}]_i > 0$, then $[{\bf b}]_i=[{\bf b}_+]_i$. 
    \item If $[{\bf b}]_i < 0$, then $[{\bf b}]_i=-[{\bf b}_{-}]_i$.
    \item $f_{\bf -b}= - f_{\bf b}$. 
    \item $f_{\bf b}^+=f_{\bf -b}^{-}$ and $f_{\bf b}^-=f_{\bf -b}^{+}$. 
    \item For any $S\subseteq \supp({\bf b})$, we have $f_{{\bf b}\vert_S}^{+}\vert f_{\bf b}^{+}$ and $f_{{\bf b}\vert_S}^{-}\vert f_{\bf b}^{-}$.
    \end{enumerate}
\end{remark}

\begin{definition} \label{sec2def1}
  For ${\bf a, b} \in \ZZ^l$, we define ${\bf a} \prec {\bf b}$, if $\supp({\bf a}_{+})\subseteq \supp({\bf b}_{+})$ and $\supp({\bf a}_{-})\subseteq \supp({\bf b}_{-})$. For binomials $f_{\bf a}, f_{\bf b}$, we define $f_{\bf a} \prec f_{\bf b}$. if ${\bf a} \prec {\bf b}$, equivalently, $\supp(f_{\bf a}^+)\subseteq \supp(f_{\bf b}^+)$ and 
  $\supp(f_{\bf a}^-)\subseteq \supp(f_{\bf b}^-)$. 
\end{definition}

\noindent 
For any $S\subseteq \{1,\ldots,l\}$ and ${\bf b}\in \ZZ^l$, we denote ${\bf b}\vert_{S}$ be the vector in $\ZZ^{\vert S\vert}$ such that $[{\bf b}\vert_{S}]_i=[{\bf b}]_i$ for all $i \in S$. For simplicity, instead of ${\bf b}\vert_{\{i:e_{i}\in E(D)\}}$, we use ${\bf b}\vert_{ E(D)}$ throughout this paper where $D$ denotes graph. 
\noindent 
Below proposition is a graph version of \cite[Proposition 4.13]{s95}. Note that $A(H)$ is a matrix obtained from $A(D)$ by deleting some columns.
\begin{proposition} \label{sec2pro1}
   Let $H$ be a oriented subgraph of a weighted oriented graph $D$ such that $V(D)=V(H)$. Then 
   \begin{enumerate}
       \item[(i)] $I_H=I_D\cap K[e_i: e_i \in E(H)]$, 
       \item[(ii)] $\mathcal{C}_H=\mathcal{C}_D\cap K[e_i: e_i \in E(H)]$,
       \item[(iii)]$\mathcal{U}_H=\mathcal{U}_D\cap K[e_i: e_i \in E(H)]$,
       \item[(iv)] $Gr_H=Gr_D\cap K[e_i: e_i \in E(H)]$. 
   \end{enumerate}
\end{proposition}

Below we recall a lemma which we use in many proofs.
\begin{lemma} \label{lem7}\cite[Lemma 3.2]{nr23}
 Let $D$ be any weighted oriented graph and $f_{\bf m} \neq 0 \in I_D$. Let $v \in V(D)$ of degree $n$. If $(n-1)$ edges of $D$ incident to $v$ are not in $\supp(f_{\bf m})$, then the other edge incident to $v$ is not in $\supp(f_{\bf m})$. Moreover if the edge $e_i$ incident to $v$ belongs to $\supp(f_{\bf m}^{+})(or \supp(f_{\bf m}^{-}))$, then there exists an edge $e_j$ incident to $v$ belongs to $\supp(f_{\bf m}^{-})(or \supp(f_{\bf m}^{+}))$.  
\end{lemma}

 \section{Strongly robustness of toric ideals of certain weighted oriented graphs} \label{sec3}

 In this section, we prove our main results. First, we exhibit a class of monomial ideals in $R$ whose toric ideals are strongly robust. We further show that the toric ideals of certain classes of weighted oriented graphs are strongly robust. Since support of a primitive binomial corresponds to connected subgraph, so here we consider connected graph.
 

    
Now we prove a lemma about the dimension of the null space of $A(D)$, where $D$ is a WOG of the following structure.

\begin{lemma}\label{sec3lem6}
Let $D$ be a weighted oriented connected graph such that $D$ has no cycles sharing a path in $D$. Let $D^{\prime}$ be the weighted oriented graph comprised of $D$ and an unbalanced cycle ${\C}$ such that ${\C}$ is connected by a new path with $D$. Let $l$ be the number of cycles in $D^{\prime}$. Then $\dim_{\QQ}{\Null}(A(D^{\prime}))=l-1$.   
\end{lemma}
\begin{proof}
Let $P_1$ be the path in $D^{\prime}$ connecting $\C$ to a vertex $v$ in $D$.
To assign vertex and edge orderings, we label the cycles and paths in the following order:  \\ 
(i)  First, label all vertices and edges of the cycle $\C$. Then label the vertices and edges of  $P_1$. \\ 
(ii)  Next, label all vertices and edges of the cycles and paths in $D$ that contain $v$ in their vertex set. Skip this step if no such cycle exists.  \\ 
(iii) Then, label all vertices and edges of cycles and paths that are connected to any vertex of the cycles labeled in (ii) via paths containing no vertices from cycles not labeled in (ii). Skip this step if no such cycle exists. \\ 
(iv) Now, label all vertices and edges of cycles and paths connected to any vertex of the cycles labeled in (iii) via paths containing no vertices from cycles not labeled in (ii) or (iii). Skip this step if no such cycle exists. \\
(v) Continue this process, labeling vertices and edges of all cycles and paths connected to those labeled in the previous step, until everything is covered.    

\underline{\bf{Claim}}: $|V(D^{\prime})|=|E(D^{\prime})|-(l-1)$.\\
\underline{{\it{Proof of the claim}}}:\; We prove this by induction on $l$. Suppose $l=1$. Then $D^{\prime}$ has only one cycle $\C$. Let $P_{2}, P_{3},\ldots,P_{n}$ be paths in  $D^{\prime}$ that are connected to $P_{1}$ where one end point of $P_{i}$ is pendant vertex for $i=2,3,\ldots,n$. Note that $|V(\C)|=|E(\C)|$ and $|V(P_i)|=|E(P_i)|+1$ for all $i$. Then  $$|V(D^{\prime})|=|V({\C})|+\sum\limits_{i=1}^{n}(|V(P_i)|-1)=|E({\C})|+\sum\limits_{i=1}^{n}|E(P_i)|=|E(D^{\prime})|.$$ 
Therefore true for $l=1$. Assume $l\geq 2$. 
Let ${\C^{\prime}}$ be one of the last step-labeling cycles (in the labeling order described above) that is connected to cycles labelled just before last step via $v^{\prime}\in V({\C}^{\prime})$. Without loss of generality, assume that paths $P_{1}^{\prime},P_{2}^{\prime},\dots,P_{n}^{\prime}$ are attached to ${\C^{\prime}}$, where one endpoint of each path is a pendant vertex. Let $D^{\prime\prime}$ be the induced subgraph of $D^{\prime}$ where $V(D^{\prime\prime})=V(D^{\prime})\setminus\left(\left(V({\C}^{\prime})\bigcup\left(\bigcup\limits_{i=1}^{n}V(P_i^{\prime})\right)\right)\setminus\{v_{i}^{\prime}\}\right)$. By induction hypothesis,\;$|V(D^{\prime\prime})|=|E(D^{\prime\prime})|-(l-2)$. Now, 
$|V(D^{\prime})|=|V(D^{\prime\prime})|+|V({\C}^{\prime})|-1+\sum\limits_{i=1}^{n}(|V(P_i^{\prime})|-1) 
=|E(D^{\prime\prime})|-(l-2)+|E({\C}^{\prime})|-1+\sum\limits_{i=1}^{n}|E(P_i^{\prime})|=|E(D^{\prime})|-(l-1)$. Thus inductive step is true and this proves our claim.
\vskip 0.2cm
\noindent 
Let $A(D^{\prime})$ be the incidence matrix of $D^{\prime}$ with respect to above vertex and edge ordering. We define $\widetilde{A(D^{\prime})}$ as the matrix obtained from $A(D^{\prime})$ by removing columns corresponding to the last labeling edge of each cycle except ${\C}$. 
Then $\widetilde{A(D^{\prime})}$ is a block triangular matrix of the  form   $\left[
\begin{array}{c|c}
  A & B \\
  \hline
  O & C
\end{array}
\right]$ where $A$ is an $|V({\C})|\times|E({\C})|$ matrix and $C$ is an upper triangular square matrix of size $|E(D^{\prime})|-(l-1)-|E({\C})|$, $O$ denotes zero matrix. Note that $A$ is the incidence matrix of ${\C}_{}$. Since the toric ideal of unbalanced cycle is zero, det($A$) $\neq$ 0. On the other hand, $C$ is an upper triangular matrix having each diagonal entry positive. This implies that det($C$) $\neq$ 0. These two arguments show that det($\widetilde{A(D^{\prime})}$) $\neq$ 0. So, rank of the matrix $\widetilde{A(D^{\prime})}$ is $|V(D^{\prime})|=|E(D^{\prime})|-(l-1)$. Therefore rank$(A(D^{\prime})) \ge |V(D^{\prime})|=|E(D^{\prime})|-(l-1)$. Again rank$(A(D^{\prime})) \le |V(D^{\prime})|=|E(D^{\prime})|-(l-1)$. This implies that rank$(A(D^{\prime})) = |V(D^{\prime})|=|E(D^{\prime})|-(l-1)$. Hence $\dim_{\QQ}$Null$(A(D^{\prime}))= l-1$.
\end{proof}
\begin{proposition} \label{pro1}
 Let $D$ be weighted oriented graph and $A(D)=[{\bf a_{1}}\; {\bf a_{2}}\;\ldots \;{\bf a_{|E(D)|}}]$ be the incidence matrix of $D$. Let $S$ be an affine semigroup in ${\ZZ_+}^{|V(D)|}$ minimally generated by $\bf a_{1},\bf a_{2},\ldots,\bf a_{|E(D)|}$. Then the Krull dimension of $K[S]$ is equal to rank of $A(D)$.
\end{proposition}
\begin{proof} Let $\psi:{\ZZ}^{|E(D)|}\rightarrow {\ZZ}^{|V(D)|}$ be a group homomorphism defined by ${\bf e_{i}} \mapsto {\bf a_{i}}$, where $\{\bf e_1,e_2,\ldots,e_{|E(D)|}\}$ is the canonical basis of ${\ZZ}^{|E(D)|}$. Then by \cite[Proposition 7.5]{ms05}, the Krull dimension of $K[S]$ is equal to $|E(D)|-$ rank(ker$(\psi)$). Since  ker$(\psi)$ is a submodule of a free module over $\ZZ$, we have ker$(\psi) \cong$ ${\ZZ}^{r}$, for some $r\le |E(D)|$. Then by Rank-Nullity theorem, we have rank$(A(D))=\vert E(D) \vert- r$, as required.
\end{proof}
\begin{corollary}
Let $D$ be a weighted oriented graph as in Lemma \ref{sec3lem6}. Then the Krull dimension of $K[e: e\in E(D)]/I_D$ is equal to $|E(D)|-(l-1)$. 
\end{corollary}
\begin{proof}
Let $A(D)=[{\bf a_{1}}\; {\bf a_{2}}\;\ldots \;{\bf a_{|E(D)|}}]$ be the incidence matrix of $D$. Let $A$ be the affine semigroup in ${\ZZ_+}^{|V(D)|}$ minimally generated by $\bf a_{1},\bf a_{2},\ldots,\bf a_{|E(D)|}$. Then $K[A]\cong \frac{K[e:e\in E(D)]}{I_D}$. Using Proposition \ref{pro1}, we have the Krull dimension of $K[A]$ is equal to the rank of $A(D)$. Then using Lemma \ref{sec3lem6}, this corollary follows.  
\end{proof}
  
\begin{lemma} \label{sec3lem5}
Let $D$ be weighted oriented graph connected to a balanced cycle ${\C_{1}}$ by a path of length $k\ge 0$ ($k=0$ means ${\C_1}$ shares only one vertex with $D$). If $f_{\bf a}\in I_D$ such that $E({\C_1})\subseteq{\supp}(f_{\bf a})$, then ${\bf a}\vert_{\{i:e_{i}\in E({\C_1})\}}\in{\Null}(A({\C_1}))$, $f_{\bf {c_1}}^{+}\vert f_{\bf a}^{+}$ and $f_{\bf{c_1}}^{-}\vert f_{\bf a}^{-}$, for some primitive binomial  $f_{\bf{c_1}}\in I_{\C_1}$ and $f_{\bf c_1} \prec f_{\bf a}$.  
\end{lemma}
\begin{proof}
Let $V(\mathcal{C}_1)=\{v_1,v_{2},\ldots,v_{m_1}\}$, $E(\mathcal{C}_1)=\{e_{1},e_{2},\ldots,e_{m_1}\}$, where $v_1$ is vertex belonging to the connecting path, and $e_{i}$ is incident to $v_{i}$ and $v_{i+1}$, for $i=1,2,\ldots,m_1$, under convention that $v_{m_1+1}=v_1$. Let $A({\mathcal{C}_1})=[a_{kl}]_{m_1\times m_1}$. Assume $f_{\bf a}\in I_D$ and  $E({\C_1})\subseteq\mbox{supp}(f_{\bf a})$. Without loss of any generality, assume that $e_{1}\in\mbox{supp}(f_{\bf a}^{+})$. 
Then we have $A(D){\bf a}=0$ which gives that for $2\le k \le m_1$,  
\begin{equation} \label{eq15}a_{k(k-1)}[{\bf a}_{+}]_{k-1}=a_{kk}[{\bf a}_{-}]_{k}\;\; \mbox{ if} \;k\; \mbox{even },
\end{equation} 
\begin{equation} \label{eq16}a_{k(k-1)}[{\bf a}_{-}]_{k-1}=a_{kk}[{\bf a}_{+}]_{k}\;\; \mbox{if}\; k\; \mbox{odd }.
\end{equation} This implies that 
\begin{eqnarray*}
a_{11}[{\bf a}_{+}]_{1} -a_{1m_1}[{\bf a}_{-}]_{m_1} &=&a_{11}\frac{a_{22}}{a_{21}}[{\bf a}_{-}]_{2} -a_{1m_1}[{\bf a}_{-}]_{m_1} (\text{by using equation \eqref{eq15} with }k=2)\\
&=&a_{11}\frac{a_{22}}{a_{21}}\frac{a_{33}}{a_{32}}[{\bf a}_{+}]_{3} -a_{1m_1}[{\bf a}_{-}]_{m_1} (\text{by using equation \eqref{eq16} with }k=3)\\
&& (\text{by repeatedly using the above step and equations \eqref{eq15},\eqref{eq16}}) \\ 
&=& \frac{a_{11}a_{22}a_{33}\ldots a_{m_1m_1}}{a_{21}a_{32}\ldots a_{m_1m_{1}-1}}[{\bf a}_{-}]_{m_1}-a_{1m_1}[{\bf a}_{-}]_{m_1}\\ 
&=& \frac{a_{11}a_{22}a_{33}\ldots a_{m_1m_1}-a_{21}a_{32}\ldots a_{m_1(m_1-1)}a_{1m_1}}{a_{21}a_{32}\ldots a_{m_1(m_1-1)}} [{\bf a}_{-}]_{m_1} \\ 
      &=& \frac{\mbox{det}(A(\mathcal{C}_1))}{a_{21}a_{32}\ldots a_{m_1(m_1-1)}} [{\bf a}_{-}]_{m_1} \\ 
      &=& 0 ~~~(\text{because det}(A(\mathcal{C}_1))=0 \text{ as $\mathcal{C}_1$ is balanced}). 
\end{eqnarray*}
 Thus the above equalities give that ${\bf a}\vert_{E(\mathcal{C}_1)} \in \text{Null}(A({\mathcal{C}_1}))$. Since Null$(A(\C_1))$ has dimension $1$, let $\{\bf c_1\}$ be a basis of Null$(A(\C_1))$, where ${\bf c_1}$ with integer entries. Thus we have $I_{\C_1}=(f_{\bf c_1})$. Then by Lemma \ref{lem7}, we get that $\supp(f_{{\bf c_1}})=\supp(f_{{\bf a}\vert_{ E(\mathcal{C}_1)}})=E(\C_1)$. Since ${\bf a}\vert_{E(\mathcal{C}_1)} \in \text{Null}(A({\mathcal{C}_1}))$, we have ${\bf a}\vert_{ E(\mathcal{C}_1)}=k {\bf c_1}$ for some integer $k\neq 0$. If $k$ is positive, then $[{\bf a}_{+}]_{i}\geq [({\bf c_1})_{+}]_{i}$ and $[{\bf a}_{-}]_{j}\geq [({\bf c_1})_{-}]_{j}$, for $e_{i}, e_{j}\in E({\C_1})$. This implies that $f_{\bf {c_1}}^{+}\vert f_{\bf a}^{+}$ and $f_{\bf {c_1}}^{-}\vert f_{\bf a}^{-}$, as required. Suppose $k$ is negative. Then we can write ${\bf a}\vert_{E(\mathcal{C}_1)}=(-k) (-{\bf c_1})$. This implies that $[{\bf a}_{+}]_{i}\geq [(-{\bf c_1})_{+}]_{i}$ and $[{\bf a}_{-}]_{j}\geq [(-{\bf c_1})_{-}]_{j}$, for $e_{i}, e_{j}\in E({\C_1})$ which yield that $f_{-\bf {c_1}}^{+}\vert f_{\bf a}^{+}$ and $f_{-\bf {c_1}}^{-}\vert f_{\bf a}^{-}$. Also, $I_{\C_1}=(f_{\bf -c_1})$ by Remark \ref{sec2rmk1}(3). This completes the proof.
\end{proof}

\begin{setup} \label{sec3nota1}
Let ${\mathcal{G}}(M)=\{{\bf x}^{{\bf a}_{i}} : 1\le i\le m\}$ be the minimal monomial generating set of a monomial ideal $M\subseteq R= K[x_{1},\ldots,x_{n}]$. Let $S=K[e_1, \ldots ,e_m]$ be the polynomial ring in the variables $e_1,\ldots,e_m$. Let $\phi$ : $K[e_1, \ldots, e_m]\rightarrow K[x_1, \ldots , x_n]$ be the $K$-algebra homomorphism defined as $\phi(e_i) ={\bf x}^{{\bf a}_{i}}$. 
\end{setup}

\begin{lemma} \label{sec3lem8}
Let $M$ be a monomial ideal as in Setup \ref{sec3nota1}. Assume $|\supp({\bf x}^{{\bf a}_{j}})|=2$ for each $j=1,2,\ldots,m$. Let $f_{\bf n}\in I_{M}$ be binomial such that $e_{i}\in\supp(f_{\bf n})$. Then for $x_{j}\in \supp(\phi(e_{i}))$, we get that $x_{j}\in \supp(\phi(e_{k}))$, $e_{k}\in \supp(f_{\bf n})$ for some $k\in\{1,2,\cdots,m\}$, $k\neq i$.    
\end{lemma}
\begin{proof}
Let $A$ be the corresponding matrix whose columns are ${\bf a}_{1},{\bf a}_{2},\cdots, {\bf a}_{m}$. Then $A{\bf n}=0$ and $[{\bf n}]_{i}\neq 0$. Suppose  $x_{j}\notin \supp(\phi(e_{k}))$ for each $k\in\{1,2,\cdots,m\}$, $k\neq i$. Then from $A{\bf n}=0$ and comparing $j^{\mbox{th}}$ entries both sides, we get $[{\bf n}]_{i}=0$ which is contradiction. This completes the proof.    
\end{proof}
\begin{lemma} \label{sec3lem7}
Let $M$ be a monomial ideal as in Setup \ref{sec3nota1}. Assume $|\supp({\bf x}^{{\bf a}_{j}})|=2$ for each $j=1,2,\ldots,m$. For a fixed $i$, suppose there is a variable $x_l \in \supp({\bf x}^{{\bf a}_{i}})$ such that there exists a monomial ${\bf x}^{{\bf a}_{j}} \not = {\bf x}^{{\bf a}_{i}} \in\mathcal{G}(M)$ with $x_l \in {\supp}({\bf x}^{{\bf a}_{i}})\cap {\supp}({\bf x}^{{\bf a}_{j}})$ and $x_l \notin {\supp}({\bf x}^{{\bf a}_{k}})$, for all $k\neq i,j$. Let $f_{\bf m}, f_{n} \in I_M$ be primitive binomials such that $f_{\bf m}^+ | f_{n}^+$.  If $e_i\in \supp(f_{\bf m}^{-})$, then $e_i\in \supp(f_{\bf n}^{-})$ and the exponent of $e_i$ in $f_{\bf m}^{-}$ is less than or equal to the exponent of $e_i$ in $f_{\bf n}^{-}$.   
\end{lemma}
\begin{proof}
    Assume that ${\bf x}^{{\bf a}_{i}}=x_{j}^{a_{j}}x_{k}^{a_{k}}$. Then $\phi(e_{i})=x_{j}^{a_{j}}x_{k}^{a_{k}}$. Let $A$ be the matrix whose columns are the exponents of generators of $M$ respectively. Then $A{\bf m}=0$. Now by the assumptions of the lemma, there is another monomial ${\bf x}^{{\bf a}_{r}}=x_{j}^{a_{j}^{\prime}}x_{l}^{a_{l}^{}}$ in $\mathcal{G}(M)$ such that $x_{j}$ belongs to exactly both the supports of $x_{j}^{a_{j}}x_{k}^{a_{k}}$ and $x_{j}^{a_{j}^{\prime}}x_{l}^{a_{l}}$. Then $\phi(e_{r})=x_{j}^{a_{j}^{\prime}}x_{l}^{a_{l}}$. Since $A{\bf m}=0$, then by equating $j^{\mbox{th}}$ entries both sides, we get $[{\bf m}]_{i}a_{j}+[{\bf m}]_{r}a_{j^{\prime}}=0$. 
 Since $e_{i}\in\mbox{supp}(f_{\bf m}^{-})$, we get $[\bf m]_i < 0$. Then by Remark \ref{sec2rmk1}(2), we get $[\bf m]_i=-[\bf m_{-}]_i$. Therefore, $[{\bf m}_{-}]_{i}a_{j}=[{\bf m}]_{r}a_{j^{\prime}}$. This implies that $[\bf m]_r > 0$, which gives that $e_{r}\in\mbox{supp}(f_{\bf m}^{+})$. Since $\mbox{supp}(f_{\bf m}^{+})\subseteq \mbox{supp}(f_{\bf n}^{+})$, we get that $e_{r}\in\mbox{supp}(f_{\bf n}^{+})$. Now, from $A{\bf n}=0$ and equating $j^{\mbox{th}}$ entries of both sides, we get $[{\bf n}]_{i}a_{j}+[{\bf n}_{}]_{r}a_{j^{\prime}}=0$. Since $e_{r}\in\mbox{supp}(f_{\bf n}^{+})$, we get $[{\bf n}]_r >0$. Then by the equality $[{\bf n}]_{i}a_{j}+[{\bf n}_{}]_{r}a_{j^{\prime}}=0$, we get that $[{\bf n}]_i < 0$. This implies that $e_{i}\in\mbox{supp}(f_{\bf n}^{-})$. As $f_{\bf m}^{+}\vert f_{\bf n}^{+}$, we get $[{\bf m}_{+}]_{r}\le [{\bf n}_{+}]_{r}$. Then by using the equalities: $[{\bf m}_{-}]_{i}a_{j}=[{\bf m}_{+}]_{r}a_{j^{\prime}}$, and $[{\bf n}_{-}]_{i}a_{j}=[{\bf n}_{+}]_{r}a_{j^{\prime}}$, we get that $[{\bf m}_{-}]_{i}\le [{\bf n}_{-}]_{i}$, as required.
\end{proof}

\noindent 
Now we give a class of monomial ideals whose toric ideals are strongly robust. 

\begin{theorem} \label{sec3thm2}
Let $M$ be a monomial ideal in $R= K[x_{1},\ldots,x_{n}]$ with ${\mathcal{G}}(M)=\{{\bf x}^{{\bf a}_{i}} : 1\le i\le m\}$ is the minimal generating set of $M$. Assume $|\supp({\bf x}^{{\bf a}_{i}})|=2$ for each $i=1,2,\ldots,m$. Suppose for each ${\bf x}^{{\bf a}_{i}}\in {\mathcal{G}}(M)$, there is a variable $x_l$ in $\supp({\bf x}^{{\bf a}_{i}})$, such that there exists a monomial ${\bf x}^{{\bf a}_{j}} \not = {\bf x}^{{\bf a}_{i}} \in\mathcal{G}(M)$ with $x_l \in {\supp}({\bf x}^{{\bf a}_{i}})\cap {\supp}({\bf x}^{{\bf a}_{j}})$ and $x_l \notin {\supp}({\bf x}^{{\bf a}_{k}})$, for all $k\neq i,j$. Then $I_{M}$ is strongly robust.  
\end{theorem}
\begin{proof}
 We show that $Gr_{M}$ is a minimal generating set of $I_{M}$. Let $f_{\bf n}\in Gr_{M}$ be a primitive binomial such that $f_{\bf n}$ is not a part of a minimal generating set of $I_{M}$. Then there exists $f_{\bf m}\in Gr_{M}\setminus\{f_{\bf n}\}$ such that $f_{\bf m}^{+}\vert f_{\bf n}^{+}$. Then by Lemma \ref{sec3lem7}, we get $e_i\in \supp(f_{\bf n}^{-})$ and the exponent of $e_i$ in $f_{\bf m}^{-}$ is less than or equal to the exponent of $e_i$ in $f_{\bf n}^{-}$. Since $e_i \in \supp(f_{\bf m}^{-})$ is arbitrary, we get $f_{\bf m}^{-}\vert f_{\bf n}^{-}$. Thus we get $f_{\bf m}^{+}\vert f_{\bf n}^{+}$ and $f_{\bf m}^{-}\vert f_{\bf n}^{-}$. This implies that $f_{\bf n}$ is not a primitive binomial. This is a contradiction. Hence $f_{\bf n}$ belongs to a set of minimal generating set of $I_M$. This implies that $Gr_{M}$ is a minimal generating set of $I_{M}$. Therefore, $I_{M}$ is strongly robust.         
\end{proof}

\begin{corollary} \label{sec3cor4}
Let $D$ be a weighted oriented graph that every edge of $D$ meets a degree $2$ vertex. Then $I_{D}$ is strongly robust. 
\end{corollary}
\begin{proof}
 Let $\mathcal{G}(I(D))$ be a minimal generating set of $I(D)$. Then $\mathcal{G}(I(D))$ satisfies the assumptions of Theorem \ref{sec3thm2}. Therefore $I_D$ is strongly robust.     
\end{proof}

\begin{corollary} \label{sec3cor3}
Let $D$ be a weighted oriented graph that is a collection of finitely many cycles that share a single vertex. Then $I_D$ is strongly robust.
\end{corollary}
\begin{proof}
For $D$, every edge meets a vertex of degree $2$. Therefore by Corollary \ref{sec3cor4}, we have $I_D$ is strongly robust.     
\end{proof}

\noindent 
Below we show a class of strongly robust toric ideals of WOGs.  

\begin{theorem}\label{sec3thm4}
Let $D$ be a weighted oriented graph with no cycles sharing a path, where each edge of every cycle meets a vertex of degree $2$. Suppose there exist cycles ${\C_1}, {\C_2}, {\C_3}$ in $D$ such that, for $i=1,2$, ${\C_i}$ is connected to ${\C_3}$ by a path $P_i$. Assume that for each edge $e\in E(P_1)\cap E(P_2)$, one of the following holds: \\
$(i)$ $e$ meets a vertex of degree $2$, or \\
$(ii)$ ${\C_3}$ shares no vertex with any other cycle in $D$ and is not connected by any path $P$ to another cycle in $D$ with $e\notin E(P)$. \\
Then $I_{D^{}}$ is strongly robust.    
\end{theorem}
\begin{proof}
By Lemma \ref{sec3lem8}, the support of a primitive binomial cannot contain any edge incident to a leaf. Since such an edge contributes nothing, we may assume without loss of generality that $D$ has no leaves.
Let $f_{\bf n}\in Gr_{D^{}}$ such that $f_{\bf n}$ is not a part of minimal generators of $I_{D^{}}$. Then there exists $f_{\bf m}\in Gr_{D^{}}\setminus \{f_{\bf n}\}$ such that $f_{\bf m}^{+}\vert f_{\bf n}^{+}$. Now we want to show that $f_{\bf m}^{-}\vert f_{\bf n}^{-}$. Let $e_{i}\in \mbox{supp}(f_{\bf m}^{-})$. Then $[\bf m]_i < 0$. \vskip 1mm 
\noindent 
{\bf CASE 1:} Suppose $e_{i}$ belongs to a cycle or $e_{i}$ belongs to a connecting path $P_i$ such that $e_{i}$ meets a degree 2 vertex. Then according to given condition, $e_{i}$ meets a degree 2 vertex say, $x_l$. Let $e_j$ be the edge in $D^{}$ that meets $e_i$ at $x_l$. Then $\phi(e_i), \phi(e_j)\in\mathcal{G}(I({D^{}}))$ such that $x_l\in \mbox{Supp}(\phi(e_i))\cap \mbox{Supp}(\phi(e_j))$. Now by applying Lemma \ref{sec3lem7}, we get that $e_i\in \supp(f_{\bf n}^{-})$ and the exponent of $e_i$ in $f_{\bf m}^{-}$ is less than or equal to the exponent of $e_i$ in $f_{\bf n}^{-}$.   

 \vskip 1mm 
\noindent 
{\bf CASE 2:} Assume $e_{i}$ belongs to a connecting path $P_i$ that avoids degree $2$ vertices. Now we prove a claim. 
\vskip 0.2cm 
\noindent 
\underline{\bf CLAIM:} For any edge $e_i$ that belongs to a connecting path $P_i$ that avoids degree $2$ vertices, and $e_i\in \supp(f_{\bf m}^{-})$. Then there exists an unbalanced cycle $\C_i$ in $D$ and a subgraph $D^{\prime}$ of $D$ such that ${\C_i}$ is connected to $D^{\prime}$ by $e_{i}$ and $E(D^{\prime})=E(D)\setminus\{E({\C_i})\cup \{e_{i}\}\}$.
\vskip 0.2cm 
\noindent 
\underline{\it Proof of the CLAIM:} Let ${\C_1}, {\C_2}$ be the cycles connected by $P_{i}$, with $e_{i}\in E(P_i)$. Suppose $e_{i}$ is incident to vertices $v_{i}$ and $v_{i}^{\prime}$. We claim that at least one of $v_{i}$ or $v_{i}^{\prime}$ lies in $V({\C_1}) \cup V({\C_2})$. Suppose not. Since $e_{i}$ avoids degree $2$ vertices, there exists cycles ${\C_3},{\C_4}$ such that ${\C_3}$ is connected to $v_{i}$ and ${\C_4}$ is connected to $v_{i}^{\prime}$. Then we get the subgraph as shown in the below figure \ref{fig1}.
\begin{figure}[h!] \centering  \includegraphics[scale=0.45]{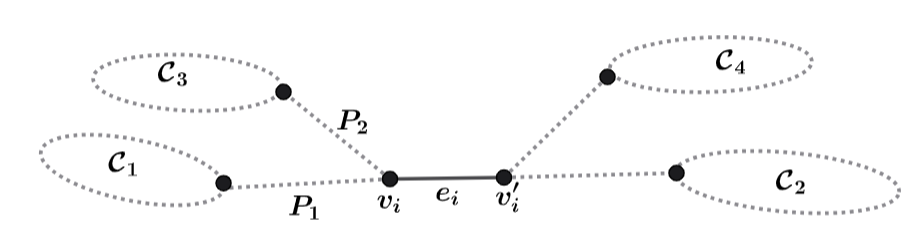}
\caption{}\label{fig1}
 \end{figure}
We see that ${\C_1}, {\C_3}$ are connected to ${\C_2}$ by paths, say $P_1, P_2$, respectively, with $e_{i}\in E(P_1)\cap E(P_2)$. Also, ${\C_2}, {\C_4}$ they are connected by a path that does not contain $e_{i}$. 
This contradicts assumption $(ii)$ of the statement. Thus at least one of $v_{i}$ or $v_{i}^{\prime}$ belongs to $V({\C_1}) \cup V({\C_2})$. Without loss of generality, assume that $v_{i}\in V({\C_1})$. We have the following two subcases.\\
\underline{\it Subcase 1:} Suppose there is no other cycle in $D$ that shares a vertex with ${\C_1}$ and is connected to ${\C_1}$ by a path not containing $e_{i}$. Then ${\C_1}$ is connected to $D^{\prime}$ via $e_{i}$, where $D^{\prime}$ is subgraph of $D$ with $E(D^{\prime})=E(D)\setminus\{E({\C_1})\cup \{e_{i}\}\}$. \\
\underline{\it Subcase 2:} Now suppose there is a cycle, say ${\C_3}$, that shares a vertex with ${\C_1}$ or is connected to ${\C_1}$ by a path not containing $e_{i}$. Then ${\C_1}, {\C_3}$ are both connected to ${\C_2}$ by paths, say $P$,$P^{\prime}$ such that $e_{i}\in E(P)\cap E(P^{\prime})$. Since $e_{i}$ does not meet a degree $2$ vertex, and by our assumption $(ii)$ there is no other cycle in $D$ such that this cycle shares a vertex with ${\C_2}$ and this cycle is connected to ${\C_2}$ by a path where $e_{i}$ not belongs to this path. Then $v_{i}^{\prime}\in V({\C_2})$. Thus ${\C_2}$ is connected to $D^{\prime}$ via $e_{i}$, where $D^{\prime}$ is the subgraph of $D$ with  $E(D^{\prime})=E(D)\setminus\{E({\C_2})\cup \{e_{i}\}\}$. 
\vskip 0.2cm 
\noindent
Thus, there exists a cycle $\C_i$ and a subgraph $D^{\prime}$ of $D$ such that ${\C_i}$ is connected to $D^{\prime}$ by $e_{i}$ and $E(D^{\prime})=E(D)\setminus\{E({\C_i})\cup \{e_{i}\}\}$. 
It remains to show that $\C_i$ is unbalanced. Since $e_{i}\in\supp(f_{\bf m}^{-})$, Lemma \ref{lem7} implies that $E({\C_i})\subseteq \supp(f_{\bf m})\cap \supp(f_{\bf n})$. We now show that ${\C_i}$ is unbalanced. Suppose, for contradiction, that  ${\C_i}$ is balanced. Then Lemma \ref{sec3lem5}, yields ${\bf m}\vert_{E({\C_i})}\in \mbox{Null}(A({\C_i}))$, so $f_{{\bf m}\vert_{ E({\C_i})}}\in I_{\C_i} \subseteq I_{D^{}}$. Note that $f_{{\bf m}\vert_{E({\C_i})}}\neq f_{\bf m}$, since $e_{i}\notin\mbox{supp}(f_{{\bf m}\vert_{E({\C_i})}})$ but  $e_i\in \text{supp}(f_{\bf m})$. This gives $f_{{\bf m}\vert_{ E({\C_i})}}^{+}\vert f_{\bf m}^{+}$ and $f_{{\bf m}\vert_{ E({\C_i})}}^{-}\vert f_{\bf m}^{-}$, contradicting the primitivity of $f_{\bf m}$. Hence, ${\C_i}$ is unbalanced. This proves the CLAIM. \\
Let $\supp(f_{\bf m})=E(D_1)$ and $\mbox{supp}(f_{\bf n})=E(D_2)$, where $D_1$ and $D_2$ are subgraphs of $D$. 

\vskip 0.2cm 
\noindent 
\underline{\bf Claim 1:} $D_1$ and $D_2$ are connected graphs. \\
{\it Proof of Claim 1:} Suppose $D_{1}$ is disconnected. Let $H_{1}$ be a connected component of $D_{1}$ such that $E(H_{1})\neq \emptyset$ and $E(D_{1})\setminus E(H_{1})\neq \emptyset$. Since $A(D^{}){\bf m}=0$, it follows that  $A(H_{1}){\bf m}\vert_{E(H_1)}=0$. Thus, $f_{{\bf m}\vert_{E(H_1)}}\in I_{H_1}$. Since $E(H_1)\neq\emptyset$ and $E(H_1)\subseteq E(D_1)=\supp(f_{\bf m})$,  it follows that  $\supp(f_{{\bf m}\vert_{E(H_1)}})\neq \emptyset$. By Lemma \ref{lem7}, $f_{{\bf m}\vert_{E(H_1)}}\neq 0$, so $f_{{\bf m}\vert_{E(H_1)}}^{+}\vert f_{\bf m}^{+}$ and $f_{{\bf m}\vert_{E(H_1)}}^{-}\vert f_{\bf m}^{-}$. This contradicts the primitivity of $f_{\bf m}$. Thus, $D_1$ is connected. This completes the proof of Claim 1. 
\vskip 0.2cm 
\noindent 
\underline{\bf Claim 2:} $\mbox{supp}(f_{\bf m})\subseteq \mbox{supp}(f_{\bf n})$. \\ 
{\it Proof of Claim 2:} It is enough to show $\mbox{supp}(f_{\bf m}^{-})\subseteq \mbox{supp}(f_{\bf n}^{-})$.  Let $e_{j}\in\mbox{supp}(f_{\bf m}^{-})$. By the given condition that each edge of every cycle in $D$ meets a degree $2$ vertex, if $e_{j}$ belongs to a cycle, then Lemma \ref{lem7} yields   $e_{j}\in\mbox{supp}(f_{\bf n}^{-})$, since $\mbox{supp}(f_{\bf m}^{+})\subseteq \mbox{supp}(f_{\bf n}^{+})$. If $e_{j}$ lies on a connecting path and meets a degree $2$ vertex, then similarly   $e_{j}\in\mbox{supp}(f_{\bf n}^{-})$. Suppose instead that $e_{j}$ lies on a connecting path but does not meet a degree $2$ vertex. By CLAIM, it follows that there is an unbalanced cycle ${\C}$ and a subgraph $D^{\prime\prime}$ such that ${\C}$ is connected by $e_{j}$ with $D^{\prime\prime}$ and $E(D^{\prime\prime})=E(D)\setminus\{E({\C})\cup \{e_{j}\}\}$. Then using Lemma \ref{lem7}, there exists an edge $e_{j^{\prime}}\in E({\C})$ and $e_{k^{\prime}}\in E(D^{\prime\prime})$ such that $e_{j^{\prime}}, e_{k^{\prime}}\in\supp(f_{\bf m}^{+})$. Then $e_{j^{\prime}}, e_{k^{\prime}}\in\supp(f_{\bf n}^{+})\subseteq E(D_2)$. Then Lemma \ref{lem7} yields $e_j\in \mbox{supp}(f_{\bf n}^{-})$. Thus $\supp(f_{\bf m})\subseteq \supp(f_{\bf n})=E(D_2)$. This completes the proof of Claim 2. 
\vskip 0.2cm
\noindent 
To complete the proof of the theorem, it suffices to show that the exponent of $e_i$ in $f_{\bf m}^{-}$ is at most the exponent of $e_i$ in $f_{\bf n}^{-}$. Let the cycles in $D_{1}$ be ${\C_1},{\C_2},\dots,{\C_{l^{\prime}}}$. 
Suppose all cycles in $D_1$ except ${\C_i}$ are balanced, with generators $f_{\bf c_{j}}$ for each ${\C_{j}}^{}$. By Lemma \ref{sec3lem6}, $\dim_{\QQ}$Null$(A(D_1))=l^{\prime}-1$ and $\{{\bf c_{j}} : j\in\{1,2,\cdots,l^{\prime}\}\setminus \{i\}\}$ forms a basis for $\Null(A(D_1))$. Since ${\bf m}$ is a $\QQ$-linear combination of these ${\bf c_{j}}$ and $e_{i}\notin E({\C}_{j}^{})$ for all such $j$, it follows that $[{\bf m}]_{i}=0$, a contradiction. Thus, at least one cycle ${\C}_{t}^{}$ in $D_{1}$, distinct from ${\C_i}$, must be unbalanced. Furthermore, $E({\C}_{t}^{})\subseteq \mbox{supp}(f_{\bf m})$. 

By Proposition \ref{sec2pro1}, $f_{\bf m}\in I_{D_2}$. Thus, the cycles $\C_i$ and $\C_t$ belong to both $D_1$ and $D_2$. Let all cycles of $D_{2}$ be ${\C^{}_1},{\C^{}_2},\dots,{\C^{}_l}$. For any $\C^{}_j\in \{\C^{}_1,\dots,\C^{}_l\}$, the toric ideal of the weighted oriented graph consisting ${\C_t^{}}$ and $\C^{}_j$ connected by a path is principal (by \cite[Theorem 5.1, Corollary 5.3]{bklo}), generated by a primitive binomial $f_{\bf a_j}$ for some vector ${\bf a_{j}}$. By Lemma \ref{sec3lem6}, we obtain $\dim_{\QQ}$Null$(A(D_2))=l-1$. Hence, $\mathcal{B}=\{{\bf a_{j}},{\bf a_{i}} : j\in\{1,2,\cdots,l\}\}$ is a basis of $\text{Null}(A(D_{2}))$. Since $E({\C_i})\subseteq \mbox{supp}(f_{\bf m})$, without loss of generality, assume that ${\bf a_i}\vert_{E({\C_i})}\prec {\bf m}\vert_{E({\C_i})}$. As ${\bf m}\vert_{E({\C_i})}\prec {\bf n}\vert_{E({\C_i})}$, it follows that ${\bf a_i}\vert_{E({\C_i})}\prec {\bf n}\vert_{E({\C_i})}$. 

Since $\mathcal{B}$ is a basis of  $\text{Null}(A(D_{2}))$ and ${{\bf m}\vert_{ E(D_2)}}, {\bf n}\vert_{E(D_2)} \in \text{Null}(A(D_{2}))$, there exist ${\mu_j},{\lambda_j}\in{\QQ}$ such that 
$${\bf m}\vert_{ E(D_2)}=\sum\limits_{j=1}^{l}\mu_{j}{\bf{a_j}}+\mu_{i}{\bf{a_i}} \text{ and } {\bf n}\vert_{ E(D_2)}=\sum\limits_{j=1}^{l}\lambda_{j}{\bf{a_j}}+\lambda_{i}{\bf{a_i}}.$$ 
Comparing the $j^{th}$ components on both sides yields $[{\bf m}]_j=\mu_i[{\bf{a_i}}]_j$ and $[{\bf n}]_j=\lambda_i[{\bf{a_i}}]_j$ for all $e_j\in E({\C_i})$. For any such  $e_{j}$, either $j\in\mbox{supp}({\bf m}_{+})\cap\mbox{supp}({\bf n}_{+})\cap \mbox{supp}(({\bf a_i})_{+})$ or $j\in\mbox{supp}({\bf m}_{-})\cap\mbox{supp}({\bf n}_{-})\cap \mbox{supp}(({\bf a_i})_{-})$, which implies $\mu_{i}>0$ and $\lambda_{i}>0$. Now compare the $i^{th}$ components: $[{\bf m}]_i=\mu_i[{\bf{a_i}}]_i$ and $[{\bf n}]_i=\lambda_i[{\bf{a_i}}]_i$. Since $\mu_{i}>0, \lambda_{i}>0$, and $e_{i}\in\mbox{supp}(f_{\bf m}^{-})$, it follows that $e_{i}\in\mbox{supp}(f_{\bf n}^{-})$. Let $v_{i}^{} \in V({\C_i})$ be the vertex incident to $e_i$, with the other two edges in $E({\C_i})$ incident to $v_{i}^{}$ denoted $e_{x}$ and $e_{y}$. As $e_i\in \text{supp}(f_{\bf m}^{-})$ and $v_{i^{}}$ is incident to the edges $e_i,e_x,e_y$, Lemma \ref{lem7}  implies that either $e_x \in$ supp$(f_{\bf m}^{+})$ or $e_y\in $ supp$(f_{\bf m}^{+})$. Without loss of generality, assume $e_x \in \text{supp}(f_{\bf m}^{+}) \subseteq \text{supp}(f_{\bf n}^{+})$. Comparing the $x^{th}$ components gives $0 < [{\bf m}]_x=\mu_i[{\bf{a_i}}]_x$ and $0 < [{\bf n}]_x=\lambda_i[{\bf{a_i}}]_x$. Since $f_{\bf m}^+ \vert f_{\bf n}^+$, we have $[{\bf m}]_x \leq [{\bf n}]_x$, so $\mu_i \leq \lambda_i$. Moreover, as $e_i \in $ supp$(f_{\bf m}^{-}) \subseteq \text{supp}(f_{\bf n}^{-})$, both $[{\bf m}]_i < 0$ and $[{\bf n}]_i <0$. Combined with $[{\bf m}]_i=\mu_i[{\bf{a_i}}]_i$ and $[{\bf n}]_i=\lambda_i[{\bf{a_i}}]_i$, and $0 < \mu_i \leq \lambda_i$, this yields $[{\bf m}]_i  \geq [{\bf n}]_i$, or equivalently, $[{\bf m}_{-}]_i\le [{\bf n}_{-}]_i$. Thus, the exponent of $e_i$ in $f_{\bf m}^{-}$ is at most the exponent of $e_i$ in $f_{\bf n}^{-}$.   
Hence, $f_{\bf m}^- \vert f_{\bf n}^-$. This contradicts the primitivity of $f_{\bf n}$.  Thus, $f_{\bf n}$ is among the minimal generators of $I_{D^{}}$, so $Gr_{D^{}}$ is a minimal generating set of $I_{D^{}}$.  It follows that $I_{D^{}}$ is strongly robust. 
\end{proof}

\begin{corollary} \label{cor9}
 Let $D$ be a weighted oriented graph consisting of cycles ${\C_1},{\C_2},\cdots,{\C_n}$ and a vertex $v$ disjoint from the vertex set of each cycle, such that for each $i\in \{1,\ldots,n\}$, the cycle ${\C_i}$ is connected to $v$ by a path $P_{i}$ with the property that $E(P_i)\cap E(P_j)=\emptyset$ whenever $i\neq j$. Then $I_{D}$ is strongly robust.    
\end{corollary}
\begin{proof}
If the cycles ${\C_1},{\C_2}$ are connected by paths, say $P_1, P_2$ to ${\C_3}$ with $e\in E(P_1)\cap E(P_2)$, then we see that any other cycle in $D$ is connected by path, say $P$ to ${\C_3}$ where $e\in E(P)$. Thus all conditions of Theorem \ref{sec3thm4} satisfied here. Then it follows.    
\end{proof}

\begin{corollary} \label{cor10}
 Let $D$ be a weighted oriented graph consisting of finitely many cycles that share only a single vertex $v$. Let $D^{\prime}$ be the weighted oriented graph obtained from $D$ by adding new cycles ${\C_1},{\C_2},\ldots,{\C_k}$, where each cycle ${\C_i}$ is connected to $D$ at $v$ by the path $P_i$  such that $E(P_i)\cap E(P_j)=\emptyset$ for all distinct $i,j\in\{1,2,\cdots,k\}$. Then $I_{D^{\prime}}$ is strongly robust.   
\end{corollary}
\begin{proof}
Like in above Corollary \ref{cor9}, it follows from Theorem \ref{sec3thm4}.    
\end{proof}

\noindent 
The Examples \ref{sec3example1}, \ref{sec3example3}, \ref{sec3example5} show that the assumptions of Theorem \ref{sec3thm4} cannot be dropped.

\begin{example} \label{sec3example1}
\begin{figure}[h!] \centering  \includegraphics[scale=0.45]{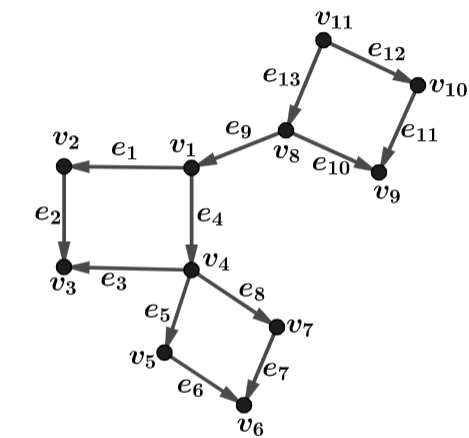}
 \caption{}\label{fig2}
 \end{figure} 
 Consider the weighted oriented graph $D$ as shown in figure \ref{fig2}. Note that the edge $e_4$ of a cycle does not meet degree $2$ vertex. Let ${\bf w}=(1,3,1,2,3,1,2,3,1,2,1)$ be the weight vector. Using Macaulay2 \cite{gs}, we get \\ 
$Gr_{D^{}}=\{e_{1}e_{3}^{3}e_{5}^{6}e_{7}^{18}e_{11}^2e_{13}-e_{2}^{3}e_{6}^{18}e_{8}^{9}e_{9}^{}e_{10}^2e_{12}, e_{1}e_{3}^{3}e_{5}^{2}e_{7}^{6}-e_{2}^{3}e_{4}e_{6}^{6}e_{8}^{3}, 
e_{4}e_{5}^{4}e_{7}^{12}e_{11}^{2}e_{13}-e_{6}^{12}e_{8}^{6}e_{9}e_{10}^{2}e_{12}, \\ 
e_{1}e_{3}^{3}e_{6}^{6}e_{8}^{3}e_{9}e_{10}^{2}e_{12}-e_{2}^{3}e_{4}^{2}e_{5}^{2}e_{7}^{6}e_{11}^{2}e_{13}, e_{1}^{2}e_{3}^{6}e_{10}^{2}e_{12}-e_{2}^{6}e_{4}^{3}e_{11}^{2}e_{13}
\}$,\\
The set of indispensable binomials $=\{e_{1}e_{3}^{3}e_{5}^{2}e_{7}^{6}-e_{2}^{3}e_{4}e_{6}^{6}e_{8}^{3},
e_{4}e_{5}^{4}e_{7}^{12}e_{11}^{2}e_{13}-e_{6}^{12}e_{8}^{6}e_{9}e_{10}^{2}e_{12}, \\ 
e_{1}e_{3}^{3}e_{6}^{6}e_{8}^{3}e_{9}e_{10}^{2}e_{12}-e_{2}^{3}e_{4}^{2}e_{5}^{2}e_{7}^{6}e_{11}^{2}e_{13}, e_{1}^{2}e_{3}^{6}e_{10}^{2}e_{12}-e_{2}^{6}e_{4}^{3}e_{11}^{2}e_{13}
\}$.
Thus $Gr_{D^{}}$ is not equal to the set of indispensable binomials i.e., $I_{D^{}}$ is not strongly robust.
\end{example}

\begin{example} \label{sec3example3}
Consider the weighted oriented graph $D$ as shown in figure \ref{fig3} consisting of cycles ${\C_1},{\C_2},{\C_3},{\C_4}$ such that each edge of each cycle in $D$ meets degree $2$ vertex but ${\C_1},{\C_2}$ are connected by paths with ${\C_3}$. 

\begin{figure}[h!] \centering  \includegraphics[scale=0.55]{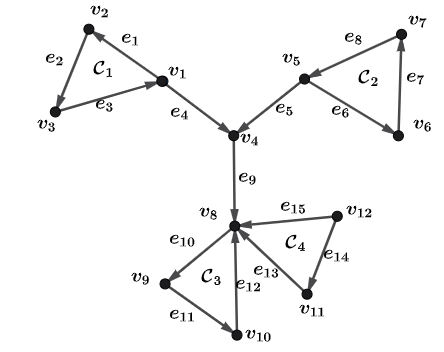}
\caption{}\label{fig3}
\end{figure}
\noindent 
 Note that the common edge $e_9$ of connecting paths does not meet degree $2$ vertex and ${\C_3}$ shares a vertex with a cycle ${\C_4}$. Let ${\bf w}=(2,1,2,3,2,2,2,2,2,2,2,1)$ be the weight vector. Using Macaulay2 \cite{gs}, we see that $Gr_{D}\neq \mbox{The set of indispensable binomials of } I_{D}$. Thus $I_{D}$ is not strongly robust.
\end{example}

\begin{figure}[h!] \centering  \includegraphics[scale=0.40]{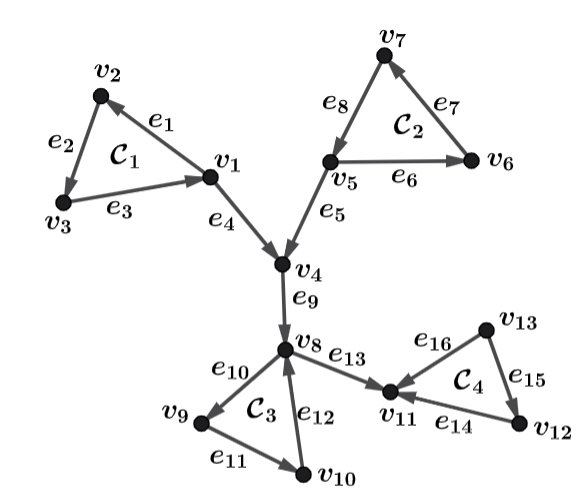}
\caption{}\label{fig4}
\end{figure} 

\begin{example} \label{sec3example5}
Consider the weighted oriented graph $D$ as shown in figure \ref{fig4} consisting of cycles ${\C_1},{\C_2},{\C_3},{\C_4}$ such that each edge of each cycle in $D$ meets degree $2$ vertex but ${\C_1},{\C_2}$ are connected by paths $P_1, P_2$ with ${\C_3}$ where the common edge $e_{9}$ of connecting paths does not meet degree $2$ vertex and ${\C_3}$ is connected by an edge $e_{13}$ with ${\C_4}$.

Let ${\bf w}=\{2,1,2,3,2,2,2,2,2,2,1,2,1\}$ be the weight vector. Using Macaulay2 \cite{gs}, we see that $Gr_{D}\neq \mbox{The set of indispensable binomials of } I_{D}$. Thus $I_{D}$ is not strongly robust.    
\end{example}

\end{document}